\documentclass{amsproc}

\usepackage{amsthm, amssymb, amsfonts}

\usepackage{graphicx}

\usepackage[dvipsnames]{xcolor} 
\usepackage[%
bookmarks=true,			
unicode=true,			
pdftitle={Deformations of noncompact manifolds},		%
pdfauthor={Elizabeth Gasparim}{Francisco Rubilar},	%
pdfkeywords={Hodge structures}{adjoint orbit}{compactification}{Lefschetz fibration},	
colorlinks=true,		
linkcolor=Black,			
citecolor=Black,		
filecolor=magenta,		
urlcolor=RoyalBlue			
]{hyperref}			
\usepackage{cleveref}

\newtheoremstyle{mydef}
{}		
{}		
{}		
{}		
{\scshape}	
{. }		
{ }		
{\thmname{#1}\thmnumber{ #2}\thmnote{ #3}}

\newtheorem{theorem}{Theorem}[section]

\newtheorem{corollary}[theorem]{Corollary}

\theoremstyle{definition}
\newtheorem{definition}[theorem]{Definition}
\newtheorem{example}[theorem]{Example}

\theoremstyle{remark}
\newtheorem{remark}[theorem]{Remark}
\newtheorem{notation}[theorem]{Notation}

\numberwithin{equation}{section}

\DeclareMathOperator{\Diag}{Diag}

\DeclareMathOperator{\Proj}{Proj}
\DeclareMathOperator{\HH}{\mathrm{H}}

\DeclareMathOperator{\Sp}{Sp}
\DeclareMathOperator{\Ad}{Ad}

\DeclareMathOperator{\Tot}{Tot}

\newcommand{\ce}{\mathrel{\mathop:}=}

\begin{document}

\title[Deformations of  Calabi--Yau manifolds and diamonds]{Deformations of noncompact Calabi--Yau manifolds, families and diamonds}


\author{Elizabeth Gasparim}
\address{Depto. Matem\'aticas, Univ. Cat\'olica del Norte, Antofagasta, Chile. }
\curraddr{}
\email{etgasparim@gmail.com}
\thanks{E. Gasparim was partially supported by Simons Associateship grant of the 
	Abdus Salam International Centre for Theoretical Physics and Network grant NT8 of the Office of External Activities, ICTP}

\author{Francisco Rubilar}
\address{Depto. Matem\'aticas, Univ. Cat\'olica del Norte, Antofagasta, Chile. }
\curraddr{}
\email{francisco.rubilar.arriagada@gmail.com}
\thanks{F. Rubilar was partially supported by Beca Doctorado Nacional Conicyt Folio 21170589 and Network grant NT8 of the Office of External Activities, ICTP}


\date{\today}

\begin{abstract}
	We discuss a new notion of deformation of complex structure, 
	which we  use as an adaptation of  Kodaira's theory of deformations,
	but that is  better suited  to the study  of noncompact  manifolds. We present several 
	families of deformations illustrating  this new approach. Our examples include 
	toric Calabi--Yau threefolds, cotangent bundles of flag manifolds, and semisimple  adjoint orbits.
	We describe their Hodge theoretical invariants, depicting Hodge diamonds and KKP diamonds.
\end{abstract}

\maketitle

\section{Deformations of complex structures}

We present a  survey of recent results on deformations  of noncompact complex manifolds 
giving a variety of examples. 
We explore some new Hodge theoretical invariants of Landau--Ginzburg models, motivated by the  KKP conjecture. 
This conjecture states that three invariants arising from  a priori unrelated fundamental ideas in complex geometry actually coincide. 
It has been recently proved that this conjecture is false in full generality, but as shown in the examples in this paper, 
there are important families where it does hold. 
Some families of deformations treated in this survey come from adjoint orbits of certain complex semisimple Lie algebras.
We focus on the noncompact case, which turns out to behave  radically different from the compact one.

Deformation theory of complex structures for compact manifolds was  initially developed by K. Kodaira and D.C. Spencer in the second half of the twentieth century. Their seminal work is summarized in the textbook ``Complex manifolds and deformations of complex structures" \cite{Ko}. We  recall some of the main results of the compact case, explain some difficulties encountered when using the same 
definition for the noncompact case,   and then state the definition we propose Def.~\ref{deformation}, which seems better suited  to the case of noncompact manifolds. 
We  give several examples initially in two and three complex dimensions, and then in further generality 
using adjoint orbits and cotangent bundles, we give examples in any complex dimension. 

\subsection{The compact case}
Let us recall some fundamental results  from the theory  of  classical 
deformations of complex structures for compact manifolds.
The intuitive idea underlying the deformation theory of complex manifolds goes as follows.
As Kodaira states in \cite[\thinspace Sec. 4.1]{Ko}, roughly speaking, a finite dimensional complex manifold  $M=M^n$ is obtained gluing domains $U_1,\ldots, U_j,\ldots$ in $\mathbb{C}^n$, i.e., $M=\bigcup_j U_j$. By \cite[Thm.\thinspace 2.1]{Ko}, these domains can be considered as polydisks. Hence, if $M$ is a compact complex manifold and $\mathfrak{U}=\{U_j\}$ is an open covering of $M$, then we can assume $\mathfrak{U}$ finite. Thus a compact complex manifold $M$ is obtained by glueing a finite number
of polydisks $U_1,\ldots,U_j,\ldots, U_l$ by identifying $z_k \in U_k$ and $z_j =f_{jk}(z_k)\in U_j: M =\bigcup_{j=1}^lU_j$ 
where $f_{jk}$ denotes the transition function associated to $U_j, U_k\in\mathfrak{U}$ such that the intersection $U_j\cap U_k$ is not empty. Hence, roughly speaking a deformation of $M$ is considered to be the glueing of the same polydisks $U_j$ via different identifications, i.e., replacing $f_{jk}(z_k)$ by $f_{jk}(z_k, t)$  where $t$ is called the parameter of deformation. Determining \textit{when} and \textit{how} the structure of $M$ as a complex manifold depends on this parameter $t$ introduced in the transition functions is the fundamental idea behind the 
concept of  deformations of complex structures.

Now we recall two basic definitions in deformation theory, namely the concepts of  family and of infinitesimal deformation.

\begin{definition}{\cite[Def. 2.8]{Ko}}\label{kodef}
	Suppose given a compact complex manifold $M_t = M_t^n$ for
	each point $t$ of a domain $B$ of $\mathbb{C}^m$.  A smooth manifold $\mathcal M = \{M_t|t\in B\}$ is called an analytic {\it family} of compact complex manifolds if there exists a holomorphic surjection $\overline{\omega}\colon \mathcal{M} \rightarrow B$ satisfying the following conditions:
	\begin{enumerate}
		\item[$\iota.$] For each $t\in B$, $\overline{\omega}^{-1}(t)$ is a compact connected subset of $\mathcal{M}$.
		\item[$\iota\iota.$] $\overline{\omega}^{-1}(t)=M_t$.
		\item[$\iota\iota\iota.$] The rank of the Jacobian matrix of $\overline{\omega}$  equals $m$ at every point of  $\mathcal{M}$.
		
		\item[$\iota\nu.$] There are a locally finite open covering $\{U_j |j = 1, 2,\ldots\}$ of $\mathcal{M}$ and complex-valued holomorphic functions $z_j^1(p),\ldots, z_j^n(p), j = 1, 2 , \ldots ,$	defined on $U_j$ such that for each $t$
		\begin{equation}
		\{p\rightarrow z_j^1(p),\ldots,z_j^n(p)|U_j\cap\overline{\omega}^{-1}(t)\neq\emptyset\}
		\label{diffamily}
		\end{equation}
		forms a system of local complex  coordinates of $M_t$.
	\end{enumerate}
\end{definition}

By $\iota.$ and  $\iota\iota.$ each $\overline{\omega}^{-1}(t)$ is a compact differentiable manifold. Condition $\iota\iota\iota.$ means that $\overline{\omega}^{-1}(t)$ is the underlying differentiable manifold of $M_t$.
With this notation $t\in B$ is the parameter of the analytic family $\{M_t|t\in B\}$ and $B$ its parameter space or base space.
This definition can be extended to the case when $B$ is an arbitrary complex manifold.
\begin{notation}
	With the same notation as before, in order to make explicit the parameter space it  is usual to write down an analytic complex family as $(\mathcal{M},\overline{\omega}, B )$.			
\end{notation}
\begin{definition}\cite[Def.\thinspace2.9]{Ko}\label{defKod} 
	Let $M$ and $N$ be two compact complex manifolds. $N$ is
	called a {\it deformation} of $M$ if $M$ and $N$ belong to the same complex analytic family, that is, if there is a complex analytic family $(\mathcal{M},\overline{\omega}, B )$ with a complex
	manifold $B$ as its parameter space such that $M =\overline{\omega}^{-1}(t_0) $ and $N = \overline{\omega}^{-1}(t_1)$
	for some $t_0,t_1 \in B$.	
\end{definition}
With the above definition in mind we get that any two elements of an analytic complex family are diffeomorphic \cite[Thm.\thinspace2.3]{Ko}. Hence, it follows that the differentiable structure of complex manifolds does not change under deformation. Several examples of deformations can be found in \cite[Sec.\thinspace 2.3]{Ko}.

The fundamental issue when we consider deformations of complex manifolds is the following: given an analytic family $(\mathcal{M},\overline{\omega}, B )$ of compact complex manifolds,
how can we determine whether the complex structure of $M_t = \overline{\omega}^{-1}$ actually
depends on $t$? Essentially the approach goes via differentiating the transition functions with respect to  the parameter $t$. Accordingly, it  is possible to show that deformations of a compact complex manifold $M$ are parametrized by 
$\HH^1(M,TM),$ that is, the first cohomology group with coefficient in its tangent sheaf. Thus, if $\theta$ is any nonzero class in $\HH^1(M,TM),$ we can associate to $\theta$ a deformation of $M$ which is called an \textit{infinitesimal deformation}. In general, computing the dimension of 
the first cohomology group with coefficients in the tangent sheaf  is not enough to prove 
existence of  deformation of a compact complex manifold. It is also necessary verify  {\it integrability}, that is, if there exist manifolds realizing such directions of deformations, and it is also necessary 
to check whether there are 
\textit{obstructions} which correspond to elements in $\HH^2(M,TM),$ \cite[Thm.\thinspace
5]{Ko}. We exhibit how this works in  some examples.

\begin{example} \label{projspace}
	The complex projective space $\mathbb{P}^n$ has no (classical) deformations of its complex structure. This result follows directly from the fact that $\HH^1(\mathbb{P}^n,T\mathbb{P}^n)=0$ and $\mathbb{P}^n$ is compact. For instance, if we consider $n=1$, we have that $T\mathbb{P}^1\cong\mathcal{O}_{\mathbb{P}^1}(2)$. Since  $\HH^1(\mathbb{P},\mathcal{O}_{\mathbb{P}^1}(2))=0$, the result follows.
\end{example}

\begin{example}
	If $M$ is a curve of genus $g$, then $\HH^2(M,TM)=0.$ Hence there are no obstructions to deforming its complex structure.
\end{example}
\begin{example}
	The Hirzebruch surface $\mathbb{F}_k$  defined as
	$$\mathbb{F}_k=\Proj(\mathcal{O}_{\mathbb{P}^1}(-k)\oplus\mathcal{O}_{\mathbb{P}^1}),\quad k\geq0$$
	is compact and has $\lfloor\frac{k}{2}\rfloor$ deformations, these are the
	$\mathbb F_{k-2i}$.
\end{example}

%

\section{Deformations of noncompact  manifolds}
In this section we describe our methods of finding infinitesimal deformations of  noncompact
manifolds which were used in \cite{BG,R} and \cite{GKRS}.
When looking for deformations of noncompact manifolds one needs to keep in
mind the caveat that cohomology calculations are generally not enough to decide
questions of existence of infinitesimal deformations, as the following example illustrates.
Suppose we  naively  extended Definition \ref{kodef}  to the noncompact case 
by simply removing the requirements of compactness, that is, repeating the same definition
only with all appearances of the word compact deleted. This would work poorly, as the following example
shows.

\begin{example} \label{ballico} Edoardo Ballico gave us the following illustration that cohomological rigidity may not imply absence of deformations,
	
	We consider deformations of $X= \mathbb C$. Clearly $\HH^1(X, TX)=0$.
	However, with the naive generalizations form the compact case, 
	there do exist nontrivial deformations of $X$ as the following family shows.
	
	Consider  $ \pi \colon  \mathbb{P}^1\times D \rightarrow D$  with $D$ any smooth manifold
	(even $\mathbb{P}^1$ or a disc) and a specific $ o\in D$.
	Take $s_\infty \colon  D \to \mathbb{P}^1 \times D$ the section of $\pi$ defined by
	$$s_\infty(x) =(\infty, x) $$
	then take another section $s$ of $\pi$ with:
	$$s(o) =(\infty,o), \qquad s(x) =(a_x,x)$$ with $a_x \in \mathbb C = \mathbb{P}^1 \setminus \{\infty\}$.
	Take as the total space for our family $Y $ as $\mathbb{P}^1 \times D$ minus the images of the two sections.
	Then we obtain a deformation of $\mathbb C$ in which at all points of $D \setminus  \{o\}$ you have $\mathbb C\setminus \{0\}$.
	Thus not a trivial deformation in any reasonable sense.
\end{example}

\begin{remark}
	Example \ref{ballico} shows that vanishing of cohomology may not imply nonexistence of deformations,
	at least not if we allow the smooth type of the manifold to vary in the family. Nevertheless, as we shall see, cohomology calculations are still useful to find deformations. 
\end{remark}In this work by deformation we mean the following:

\begin{definition}\label{deformation}
	A {\it deformation family} of a complex manifold $X$ is a holomorphic surjective submersion $\tilde{X}\stackrel{\pi}{\rightarrow}D$, where $D$ is a complex disc centred at $0$ (possibly a vector space, possibly infinite dimensional), satisfying:
	\begin{itemize}
		\item $\pi^{-1}(0)=X$,
		\item $\tilde{X}$ is locally trivial  in the $C^{\infty}$ category.
	\end{itemize}
	The fibres $X_t=\pi^{-1}(t)$ are called \textit{deformations} of $X$.
\end{definition}

In further   generality, allowing the parameter space to be a variety or a scheme, and  requiring that the bundle be locally trivial,
we obtain the concept of a {\it family} of complex structures.

\begin{remark}
	We say that a deformation is trivial if it is isomorphic in the holomorphic category to the original manifold.
\end{remark}

\begin{remark}
	Our choice for the dimension of $D$ is $n=h^1(X,TX)$ whenever possible. The case $n=0$ may be included in the following definition.
\end{remark}

\begin{definition}\label{formal}
	We call a manifold $X$ {\it formally rigid} when $\HH^1(X, TX)=0$.
\end{definition}

\begin{definition}\label{rigid}
	We call a manifold $X$ {\it rigid} if any deformation
	$\tilde{X}\stackrel{\pi}{\rightarrow}D$
	is biholomorphic to the trivial bundle $X \times D \to D$.
\end{definition}

\begin{example} Projective spaces are rigid. Indeed, Ex.~\ref{projspace} shows it is formally rigid, and in the 
	compact case formally rigid implies rigid. 
\end{example}

In general, formally rigid does not imply rigid, as example \ref{ballico} shows, but we improve on the situation 
with our proposed definition.
With Def.~\ref{deformation} we do not claim to solve the problem that a manifold $X$ does not deform under the condition $\HH^1(X,TX)=0$, however we certainly eliminate some unwanted 
pathological cases such as the one in Example \ref{ballico}.

Observe that the deformations considered in \cite{BG} satisfy Definition \ref{deformation}, hence maintain the $C^\infty$ type of the manifold. 


\subsection{Surfaces}

We give examples of noncompact surfaces which have finite dimensional spaces of deformations,
and  exhibit their  families of  deformations. Let $Z_k$ be the noncompact surface defined as
$$Z_k=\Tot(\mathcal{O}_{\mathbb{P}^1}(-k));\qquad k\geq1.$$
Barmeier and Gasparim \cite{BG}  described the classical deformation theory of the surfaces $Z_k$, they have also described their
noncommutative deformation theory in \cite{BG2}, but in this work we will restrict ourselves to classical deformations.
Here we recall a definition taken from \cite{M}, in order to clarify the exposition of the results for the non-compact surfaces $Z_k$.

These surfaces are total spaces of negative line bundles on the projective line. A preliminary
estimate  via cohomology  calculations \cite[Thm.\thinspace 5.4]{BG} showed that $Z_k$ admits a $(k-1)$-dimensional semiuniversal family of classical deformations. Here by semiuniversal  it is meant that the family is versal and the 
Kodaira--Spencer map is bijective, see \cite{M}.
Denoting by $\mathcal{Z}_k$ any nontrivial deformation of $Z_k$ for $k\geq 2$,  \cite[Thm.\thinspace 6.6]{BG} showed that $\mathcal{Z}_k$ contains  no compact complex analytic curves. Furthermore, \cite[Thm.\thinspace6.14]{BG} showed that any holomorphic vector
bundle on $\mathcal{Z}_k$ splits as a direct sum of algebraic line bundles. This is somewhat surprising, given the
existence of nontrivial moduli of vector bundles on the original $Z_k$ surfaces as constructed in \cite{BGK}. Also in \cite{BG} it is proved that any nontrivial deformation of $Z_k$ is affine (this situation contrasts with the case of the threefolds considered in the next section).

Deformation theory effects on moduli of vector bundles have applications to mathematical physics. The
effects
for the case of surfaces $Z_k$ can be interpreted in Yang--Mills theory. 
They  imply that moduli of $SU(2)$ instantons on noncompact surfaces are sensitive to the complex
structure: the moduli of irreducible instantons  of normalized charge   (or splitting) $j$
over the noncompact surfaces $Z_k$ are of dimension $2j-k-2$,
whereas a nontrivial (commutative) deformation of ${Z}_k$ admits no instantons \cite[Thm.\thinspace7.4]{BG}.
In comparison,  in the context of deformations of compact curves and surfaces and their moduli we note that by
Grothendieck's splitting theorem  any holomorphic vector bundle on $\mathbb{P}^1$ splits as a direct sum of
(algebraic) line bundles. Neither the curve $\mathbb{P}^1$ itself nor its moduli spaces of vector bundles admit any
deformations.
Curves of higher genus do admit deformations and a celebrated theorem of Narasimhan and Ramanan
\cite{NR, BV} shows that all classical deformations of the moduli of stable bundles on a smooth curve come from  deformations of the curve itself (case $g > 1,\quad (r, d) = 1$).
As we saw,  deformations of the surfaces $Z_k$  does not produce new moduli of bundles. 
The situation for Calabi--Yau threefolds turns out a lot more interesting, as we shall now see.

\subsection{Calabi--Yau threefolds}
We define  Calabi--Yau threefolds $W_k$ as total spaces of holomorphic vector bundles on $\mathbb{P}^1$, i.e.,
$$W_k=\Tot(\mathcal{O}(-k)\oplus\mathcal{O}(k-2));\qquad k\geq1.$$
We put \textit{canonical coordinates} on $W_k$ by taking the open covering $\mathfrak{U}=\{U, V\}$, where $U=\{(z, u_1,u_2)\}\simeq\mathbb{C}^3\simeq V=\{(\xi, v_1, v_2)\}$ such that  
$$(\xi, u_1, u_2)=(z,z^ku_1,z^{2-k}u_2),$$
on the intersection $U\cap V\simeq \mathbb{C}^{\ast}\times \mathbb{C}^2$.
The behaviour of the deformations of these manifolds is quite different from the one of deformations of the surfaces $Z_k$. This situation is described by the following results proved by Rubilar in \cite{R} and  Gasparim and Suzuki in \cite{GS}. First, Rubilar proved that while $W_1$ is formally rigid \cite[Lem.\thinspace 5.2.1]{R},  $W_2$ has an infinite dimensional family of deformations \cite[Thm.\thinspace 6.3.2]{R}. The latter was accomplished by showing that $\dim \HH^1(W_2,TW_2) =\infty $ and then proving that directions of deformations parametrized
by such cohomology are integrable, by explicitly constructing  the corresponding families. 
Moreover, he also showed that $W_2$ has both affine and non-affine deformations, situation which 
contrasts with the one of the surfaces $Z_k.$
Subsequently, Gasparim and Suzuki showed explicitly that there are infinitely many isomorphism types among the deformations of $W_2$  \cite[Thm.\thinspace 3.7]{GS}. 
Then, presenting  transformations taking deformations of $W_k$ for $k>2$ to deformations of $W_2$
and conversely, they obtain also infinitely many isomorphism types of deformations of $W_k$ for 
$k>2.$
\subsubsection{Moduli of vector bundles}
The study of deformations is very closely linked to the theory  of moduli spaces. 
In fact, there is a precise sense in which infinitesimal deformations can be 
used to generate neighborhoods of moduli spaces. For the case of vector bundles 
this is known as Kuranishi theory. Here we only comment on the results for our 
surfaces and threefolds, without discussing the general theory. 

For deformations of  surfaces $Z_k$, Barmeier and Gasparim  showed that moduli spaces of holomorphic vector bundles  with fixed topological invariants are trivial \cite[Cor.\thinspace 6.19]{BG},
i.e. consist of a single point. For the threefold $W_2$, Gasparim and Suzuki showed that 
some nontrivial deformations of $W_2$ have the effect of decreasing the dimension of the moduli of vector bundles  by just one dimension, hence keeping nontrivial moduli \cite[Thm. \thinspace 8.13]{GS}.

In conclusion,  we have observed
that the theories of deformations of surfaces $Z_k$ and threefolds $W_k$ 
are qualitatively different, the latter being very rich in terms of applications to moduli of vector bundles. 

\section{Hyperk\"ahler families}
\begin{definition}
	A \textit{hyperk\"ahler manifold} is a Riemannian manifold of real
	dimension $4k$ and holonomy group contained in $\Sp(k)$.
\end{definition} 

Here $\Sp(k)$ denotes a compact form of a symplectic group, identified with the group of quaternionic-linear unitary endomorphisms of a  $k$-dimensional quaternionic Hermitian space.
Every hyperk\"ahler manifold $M$ has a 2-sphere of complex structures  with respect to which the metric is K\"ahler.

In particular, it is a hypercomplex manifold, meaning that there are three distinct complex structures, $I, J,$ and $K,$ which satisfy the quaternion relations
$$\displaystyle I^{2}=J^{2}=K^{2}=IJK=-1.$$
Any linear combination
$$ aI+bJ+cK\,$$
with $ a,b,c $ real numbers such that
$ a^{2}+b^{2}+c^{2}=1$
is also a complex structure on $M$.
Hence, hyperk\"ahler manifolds are especially rich examples where to study deformations of complex structures. 
Each of the complex structures $ aI+bJ+cK$  may be regarded as a deformation of $I$.

\begin{remark}
	To start with $I, J, K$ are almost complex structures, these are then required to be integral,
	and in such case the complex structures obtained agree with the ones defined by Kodaira using charts.
\end{remark}

\begin{remark} We have the following general facts:
	\begin{enumerate}
		\item Hyperk\"ahler manifolds are special classes of K\"ahler manifolds. 
		\item  All hyperk\"ahler manifolds are Ricci-flat and are thus Calabi--Yau manifolds. 
	\end{enumerate}
\end{remark}
Compact hyperk\"ahler manifolds have been extensively studied using techniques from algebraic geometry, 
sometimes under  the alternative  name of holomorphically symplectic manifolds. Due to Fedor Bogomolov's decomposition 
theorem (1974,\cite{Bog}), the holonomy group of a compact holomorphically symplectic manifold $M$ is exactly $\mathrm{Sp}(k)$ if and only if $M$ is simply connected and any pair of holomorphic symplectic forms on $M$ are scalar multiples of each other. 

For applications to the study of Hodge theory of Landau--Ginzburg  (LG) models (see Def. \ref{defLGmodel}), here we wish to focus on those hyperk\"ahler families that contain the adjoint orbits of semisimple 
Lie groups. An adjoint orbit of a compact Lie group is called a {\it flag manifold}. 

Families of hyperk\"ahler structures containing  semisimple adjoint orbits   were studied by 
Kronheimer \cite{Kro1}, Biquard \cite{Biq}, and Kovalev \cite{Kov}. Such families include 
both cotangent bundles of flag manifolds and adjoint orbits of semisimple Lie algebras. 
In the following section, we will add superpotentials to these families and then study  them from the point of view of 
Landau--Ginzburg models. The approach we propose here goes in a different way, we study the hyperk\"ahlerian structures 
using Lie theory techniques.
We first recall  some definitions 
and results.

\begin{theorem}\cite[Thm.\thinspace1.1]{Kov}\label{teo1.1Kovalev}
	Let $G$ be a compact semisimple Lie group with $\mathfrak g$  its Lie
	algebra, $H$ a stabilizer of some non-zero element of $\mathfrak{g}$ under the adjoint
	action of $G$, and $T$ a maximal torus such that $G \supseteq H \supseteq T$. Let $G^c$, $H^c$ be
	the complex forms of $G$ and $H$.
	Then the manifold $G^c/H^c$ has a family of hyperk\"ahlerian structures
	with the following properties:
	\begin{enumerate}
		\item 	\begin{enumerate}
			\item the hyperk\"ahlerian metrics are complete,
			\item  the family admits parametrisation as follows: let $\mathfrak t$ be a Cartan
			subalgebra of $\mathfrak g$ corresponding to $T$, then the parameter space consists of
			those triples $(\tau_1, \tau_2, \tau_3)$ in $\mathfrak{t}$ which each have the stabilizer $H$;
		\end{enumerate} 
		\item  On $G^c/H^c$ there is a complex structure $I$ which is independent of $\tau_1$.
		If also the complex group $H^c$ is the stabilizer of the element $\tau^c = \tau_2 + i\tau_3$
		then the underlying complex manifold determined by $I$ is isomorphic to the
		complex adjoint orbit $\mathrm{Ad}(G^c)\tau^c$.	
	\end{enumerate}
\end{theorem}

\begin{remark}
	Observe that putting $H^c=T^c$ we obtain the result of Kronheimer \cite{Kro1} on regular semisimples orbits, that is, those of type 
	$G^c/T^c$. In general $T^c$ is  a proper subgroup of $H^c$, therefore in such case the family of hyper-K\"ahlerian structures is smaller. 
\end{remark}

To motivate the parametrisation, consider  the adjoint action of $G^c$. Let $\tau^c$ be some element of $\mathfrak{g}^c$ and $H^c$ its stabilizer. The orbit $\mathrm{Ad}(G^c)\tau^c$ is diffeomorphic to $G^c/H^c$. The second de Rham cohomology of $G^c/H^c$ can be identified with a real commutative Lie algebra $\mathfrak{h}_\tau\subset\mathfrak{g}$ determined (up to isomorphism) by \textit{regularity properties} of $\tau^c$: for example, $\mathfrak{h}_\tau=\mathfrak{t}$ if $\tau^c$ is a regular element of $\mathfrak{t}^c$. The
three K\"ahler structures on a complex orbit and the cohomology classes of their
K\"ahler forms can be put in correspondence to elements $\tau_1, \tau_2, \tau_3$ in $\mathfrak h_{\tau}$ such that $\tau_2+i\tau_3=\tau^c$. Then $\tau_1$ serves as a parameter in the family of hyperk\"ahlerian structures induced on $G^c/H^c$ from the orbit of $\tau_2+i\tau_3$.

\begin{theorem}\cite[Thm.\thinspace1.2]{Kov}\label{teo2Kov}
	Let $G$ and $\mathfrak{g}$ be as in Thm.~\ref{teo1.1Kovalev}. If $\tau_2, \tau_3\in\mathfrak{t}$ and $\sigma^c$ is a nilpotent element of the centralizer $z(\tau_2+i\tau_3)$ in $\mathfrak{g}^c$ then the complex adjoint orbit $\mathrm{Ad}(G^c)(\tau_2+i\tau_3+\sigma^c)$ has a family of hyperk\"ahlerian structures. The family is parametrised by elements $\tau_1\in\mathfrak{t}$ such that the centralizers in $\mathfrak{g}$ of the pair $(\tau_2, \tau_3)$ and the triple $(\tau_1, \tau_2, \tau_3)$ coincide.
\end{theorem}
\begin{remark}
	If $\tau_2=\tau_3=0$ then $\tau_1=0$ is the only possible value and we obtain the following result due to  Kronheimer.
\end{remark}

\begin{corollary}\thinspace\cite{Kro2}
	A nilpotent adjoint orbit of $G^c$ is a hyperk\"ahler manifold. 
	
\end{corollary}

It remains to point out that Theorem \ref{teo2Kov} applies to all  adjoint orbits of $G^c$. This follows immediately from standard properties of complex semisimple Lie algebras, the details can be found e.g. in Humphreys \cite{Hum} or San Martin \cite{SM}. Any element of $\mathfrak{g}^c$ can be written uniquely in the form $x=s+n$, with $n$ nilpotent, $s$ semisimple and $[s, n]=0$. As a semisimple element $s$ is contained in some Cartan subalgebra, that is in $\Ad(g)\mathfrak{t}^c$ for some $g\in G^c$. So any element of $\mathfrak{g}^c$ is $G^c$-conjugate to one of the form $(\tau_2+i\tau_3)+\sigma^c$.


We state in a simplified form of a consequence of Thm.~\ref{teo1.1Kovalev} which we will use to study deformations.

\begin{corollary}\label{deforb} Let $\mathfrak g$ be a semisimple complex Lie algebra, and $H_o\in \mathfrak g $ a regular element.
	Let $G$ be a connected Lie group with Lie algebra $\mathfrak g$ and $K$ a compact form of $G$. 
	Then the adjoint orbit $\Ad G (H_0) $ is a deformation of the cotangent bundle of the flag manifold $\Ad K (H_0).$
\end{corollary}

\section{Adjoint orbits and cotangent bundles}

Using the result of Thm.~\ref{teo1.1Kovalev} we see that adjoint orbits may be regarded as deformations of cotangent bundles of flag manifolds. 
That is, if $G$ is a complex Lie group and $K$ its compact part, then an adjoint orbit $\Ad G (H_o)$ can be regarded as 
a deformation of the cotangent bundle of the flag manifold $\Ad K(H_0)$. 
We wish to add   superpotentials to these families,  to look at them  from the point of view of Landau--Ginzburg models, 
and then to study  how their Hodge diamonds vary in the families. We first recall the construction  
given in \cite{GGSM1}.

\begin{definition} \label{defLGmodel}\emph{A Landau--Ginzburg model}  (LG) is a pair $(Y,w)$, where
	
	\begin{enumerate}
		\item $Y$ is a smooth complex quasi-projective variety with trivial canonical bundle $K_Y$;
		
		\item $w\colon Y\to \mathbb C$ is a morphism with a compact critical locus $\mathrm{crit}(w)\subset Y$, 
		called the superpotential. 
	\end{enumerate}
\end{definition}

\subsection{LG models on adjoint orbits}

Let $\mathfrak g$ be a complex semisimple Lie algebra with Cartan subalgebra $\mathfrak h$, and $\mathfrak h_\mathbb R$
the real subspace generated by the roots of $\mathfrak h$. Let $\Pi$ denote the set of all roots of $\mathfrak{h}$.
An element $H\in \mathfrak{h}$ is  called \textit{regular} if $\alpha \left(
H\right) \neq 0$ for all $\alpha \in \Pi $.

\begin{notation} $\mathcal{O}\left(H_{0}\right)\ce \Ad G (H_o)$.
\end{notation}

\begin{theorem}\label{thm1}
	\cite[Thm.\thinspace3.1]{GGSM1}
	Given $H_{0}\in \mathfrak{h}$ and $H\in \mathfrak{h}_{\mathbb{R}}$ with $H$ a regular  element, the potential $f_{H}:\mathcal{O} \left( H_{0}\right) \rightarrow \mathbb{C}$ defined by
	\[
	f_{H}\left( x\right) = \langle H,x\rangle \qquad x\in \mathcal{O}\left(H_{0}\right)
	\]
	has a finite number of isolated singularities and defines a symplectic Lefschetz fibration; that is to say
	\begin{enumerate}
		\item the singularities are (Hessian) nondegenerate;
		\item if $c_{1},c_{2}\in \mathbb{C}$ are regular values then the level manifolds $f_{H}^{-1}\left( c_{1}\right) $ and $f_{H}^{-1}\left( c_{2}\right) $ are diffeomorphic;
		\item there exists a symplectic form $\Omega $ on $\mathcal{O}\left(H_{0}\right) $ such that the regular fibres are symplectic submanifolds;
		\item each critical fibre can be written as the disjoint union of affine subspaces contained in $\mathcal O \left( H_0 \right)$, each symplectic with respect to $\Omega$.
	\end{enumerate}
\end{theorem}

Given the Iwasawa decomposition: $\mathfrak g = \mathfrak k \oplus  \mathfrak a \oplus  \mathfrak n$, 
if $\Pi$  is a set of roots of $\mathfrak a$, with a choice of a set of positive roots $\Pi^+$ and simple roots
$\Sigma \in \Pi^+$, then the 
corresponding Weyl chamber is $\mathfrak a^+$.
A subset $\Theta \in \Sigma$  defines a parabolic subalgebra $\mathfrak p_\Theta$  with parabolic subgroup 
$P_\Theta$  and a flag $ \mathbb{F}_{\Theta }=G/P_\Theta$.
A element 
$H_\Theta  \in  cl \mathfrak a^+$ is characteristic for $\Theta \subset \Sigma $ if $\Theta  = \{\alpha \in  \Sigma : \alpha (H_\Theta) = 0\}$. 

\begin{theorem}\cite[Thm.\thinspace 2.1]{GGSM2}\label{iso}
	\label{teodifeocotan}The adjoint orbit $\mathcal{O}\left( H_{\Theta }\right)
	=\mathrm{Ad}\left( G\right) \cdot H_{\Theta }\simeq G/Z_{\Theta }$ of the
	characteristic element $H_{\Theta }$ is a $C^{\infty }$ vector bundle over $%
	\mathbb{F}_{\Theta }$ isomorphic as a vector bundle to the cotangent bundle $T^{\ast }\mathbb{F}%
	_{\Theta }$. Moreover, we can write down a diffeomorphism $\iota :\mathrm{Ad%
	}\left( G\right)\cdot H_{\Theta }\rightarrow T^{\ast }\mathbb{F}_{\Theta }$
	such that
	
	\begin{enumerate}
		\item $\iota $ is equivariant with respect to the actions of   $K$, that is,
		for all $k\in K$,
		$$\iota \circ \mathrm{Ad}\left( k\right) =\widetilde{k}\circ \iota$$
		where $K$ is the compact subgroup in the Iwasawa decomposition $G=KAN$, and
		$\widetilde{k}$ is the lifting to $T^{\ast }\mathbb{F}_{\Theta }$ (via
		the differential) of the action of $k$ on $\mathbb{F}_{\Theta }$.
		
		\item The pullback of the canonical symplectic form on $T^{\ast }\mathbb{F}%
		_{\Theta }$ by $\iota $ is the (real) Kirillov--Kostant--Souriau form on
		the orbit.
	\end{enumerate}
\end{theorem}

Viewing the orbit as the cotangent bundle of a flag manifold, we can identify the topology  of the fibres in terms of the topology of the flag.
Denote by $\mathcal{W}$ the Weyl group of $G$, then we have the following two corollaries.
\begin{corollary}\cite[Cor.\thinspace 4.5]{GGSM1}
	The homology of a  regular fibre  coincides  with the homology of
	$\mathbb{F}_{\Theta }\setminus \mathcal{W}\cdot H_{\Theta}$.
	In particular the middle Betti number is $k-1$ where $k$ is
	the number of singularities of the fibration (equal to the number of elements in
	$\mathcal W \cdot H_\Theta$).
\end{corollary}

For the case where singular fibres have only one critical point, we have the following corollary.

\begin{corollary}\cite[Cor.\thinspace 5.1]{GGSM1}
	\label{cor.sing}
	The homology of the singular fibre though $ w H_\Theta$, $w \in \mathcal{W}$,
	coincides with that of
	\[
	\mathbb{F}_{H_\Theta} \setminus \{uH_\Theta \in \mathcal{W}\cdot H_{\Theta} | u \neq w\}.
	\]
	In particular, the middle  Betti number of this singular fibre
	equals $k-2$, where $k$ is the number of singularities of the  fibration $f_H$.
\end{corollary}

These corollaries show that Hodge diamonds for the LG models can have arbitrarily high numbers in 
their middle cohomology, and that these Lefschetz fibrations may have large quantities of vanishing cycles.

\section{Diamonds}

\subsection{Hodge diamonds}
We now wish to consider some examples of how diamonds vary under deformations. 
Consider first the standard Hodge diamond of a variety:
\begin{displaymath}
\begin{array}{cccccccccc}
& & & & h^{n,n} & & & & & \\
& & & & & & & & \\
& & & h^{n,n-1} & & h^{n-1,n} & & & & \\
& & & & & & & & & \\
& & h^{n,n-2} & & h^{n-1,n-1} & & {h^{n-2,n}} & & & \\
& \reflectbox{$\ddots$} & & & \vdots & & & {\ddots} & & \\
h^{n,0} & & & {\cdots} & {\stackrel{\curvearrowleft}{{}}}
& {\cdots} & & & { h^{0,n}} & \updownarrow { }\star\\
& \ddots & & & {\vdots} & & &\reflectbox{$\ddots$} & & \\
& & h^{2,0} & & {h^{1,1}} & & {h^{0,2}} & & & \\
& & & & & & & & & \\
& & & {h^{1,0}} & & {h^{0,1}} & & & & \\
& & & & & & & & & \\
& & & & {h^{0,0}} & & & & &\\
& & & & \stackrel{\longleftrightarrow}{{}} & & & & &
\end{array}
\end{displaymath}
where the symbols $\curvearrowleft, \longleftrightarrow, \updownarrow { }\star$ 
are there to remind us that in the case of  smooth projective varieties
there are symmetries of the diamond corresponding to 
Serre duality, conjugation, and Hodge star, respectively. 

\subsection{KKP diamonds}

For the case of Landau--Ginzburg models $(Y,w)$ there are three new Hodge theoretical invariants which were defined by 
Katzarkov, Kontsevich, and Pantev in \cite{KKP}. These invariants take into account 
not just the variety, but also the potential together with its 
critical points and vanishing cycles. They are the numbers $f^{p,q}(Y, w)$  which come from
sheaf cohomology of  logarithmic forms, the
numbers $h^{p,q}(Y, w)$  motivated by mirror symmetry considerations, and the numbers 
$i^{p,q}(Y, w)$  defined using ordinary mixed Hodge theory.
%
%
%
%

To define these new invariants, the authors require that the LG model be tamely compactifiable in 
the following sense.

\begin{definition}\cite{KKP} \label{def-3}\emph{A tame compactified Landau--Ginzburg model} is the data $((Z,f),D_Z)$, where
	
	\begin{enumerate}
		\item $Z$ is a smooth projective variety and $f\colon Z\to \mathbb P ^1$ is a flat morphism.
		
		\item $D_Z=(\cup _i D^h_i)\cup (\cup _jD_j^v)$ is a reduced normal crossings divisor such that
		
		\begin{itemize}
			\item[(i)] $D^v=\cup _jD^v_j$ is a scheme theoretical pole divisor of $f$, i.e. $f^{-1}(\infty)=D^v$. In particular $ord _{D^v_j}(f)=-1$ for all $j$;
			
			\item[(ii)] each component $D_i^h$ of $D^h=\cup _iD^h_i$ is smooth and horizontal for $f$, i.e. $f\vert _{D^h_i}$ is a flat morphism;
			
			\item[(iii)] The critical locus $crit(f)\subset Z$ does not intersect $D^h$.
			
		\end{itemize}

		\item $D_Z$ is an anticanonical divisor on $Z$.
		
		\noindent One says that $((Z,f),D_Z)$ is \emph{a compactification of the Landau--Ginzburg model}
		$(Y,w)$ if in addition the following holds:
		
		\item $Y=Z\setminus D_Z$, $f\vert _Y=w$. 
	\end{enumerate}
	
\end{definition}


\begin{remark} A caveat about algebraic compactifications should be mentioned here. 
	\cite{BCG} showed that the choice of compactification may have strong effects 
	on the Hodge diamonds.  In fact, they give examples when  
	the topology of the  compactified regular fibre for $f_H$ changes drastically 
	according to the choice of homogenisation of the ideal cutting out the orbit as an affine variety. 
	For the case of the maximal adjoint orbit of $\mathfrak{sl}(3, \mathbb C)$, namely the one 
	diffeomorphic to the cotangent bundle of the full flag $F(1,2)$, \cite{BCG} give examples 
	of two algebraic compactifications of such an orbit which produce in one case 
	$h^{1,4}=h^{4,1}=16$ and in the other case $h^{1,4}=h^{4,1}=1$. Such radical difference 
	being produced simply by the choice of homogenisation of the ideal defining the 
	adjoint orbit. Thus, one must be very careful when using compactifications to study the 
	Hodge theory of noncompact varieties. This remark also highlights the importance of the 
	careful definition of a tame compactification.
\end{remark}

Assume that we are given a Landau--Ginzburg model $(Y,w)$ with a tame compactification~$((Z,f),D_Z)$ as above.
Denote by $n=\dim Y=\dim Z$ the (complex) dimension of~$Y$ and $Z$. Choose a point $b\in  \mathbb C$ which is near $\infty$ 
and such that the fibre $Y_b=w^{-1}(b)\subset Y$ is smooth.
Let us briefly recall the definitions, for more details see \cite{LP}.

\subsubsection{$f^{p,q}(Y,w)$}\label{subs-def-fpq} 
Let $\Omega^\bullet _Z(log\, D_Z)$
denote  the logarithmic de Rham complex.
For each $a\geq 0$ define the  \emph{a sheaf $\Omega ^a _Z(log\,D_Z ,f)$ of $f$-adapted logarithmic forms} as:
$$\Omega ^a _Z(log\,D_Z ,f)=\{\alpha \in \Omega ^a _Z(log\,D_Z)\ \vert \ df\wedge \alpha \in \Omega ^{a+1} _Z(log\,D_Z )\}.$$

\begin{definition} \label{fpq}\emph{The Landau--Ginzburg Hodge numbers} $f^{p,q}(Y,w)$ are defined as follows:
	$$f^{p,q}(Y,w)=\dim \HH^p(Z,\Omega ^q _Z(log\,D_Z ,f)).$$
\end{definition}

\subsubsection{$h^{p,q}(Y,w)$} 
\begin{definition}
	\label{hpq} \cite[Def.~8]{LP} Assume that $(Y,w)$ is a Landau--Ginzburg model of Fano type. Consider the relative cohomology $\HH^\bullet(Y,Y_b)$ with the nilpotent operator $N$ and the induced canonical filtration $W$. \emph{The Landau--Ginzburg numbers} $h^{p,q}(Y,w)$ are defined as follows:
	$$h^{p,n-q}(Y,w)=\dim gr _{2(n-p)}^{W,n-a}\HH^{n+p-q}(Y,Y_b)\ \ \text{if $a=p-q\geq 0$},$$
	$$h^{p,n-q}(Y,w)=\dim gr _{2(n-q)}^{W,n+a}\HH^{n+p-q}(Y,Y_b)\ \ \text{if $a=p-q< 0$}$$
	where $gr$ is the standard grading of the mixed Hodge structure by the weight filtration. 
\end{definition}

\subsubsection{$i^{p,q}(Y,w)$} For each $\lambda \in \mathbb C$ one
has the corresponding sheaf~$\phi _{w-\lambda}\mathbb C _Y$ of vanishing cycles for the fibre $Y_\lambda$.
The sheaf $\phi _{w-\lambda}\mathbb C _Y$ is supported on the fibre $Y_\lambda$ and is equal to zero if $\lambda $ 
is not a critical value of $w$. 

\begin{definition}\label{ipq}
	Assume that the horizontal divisor $D^h\subset Z$ is empty, i.e. assume that the map $w\colon Y\to \mathbb C$ is proper. Then
	\begin{enumerate}
		\item
		\emph{the Landau--Ginzburg Hodge numbers}~$i^{p,q}(Y,w)$ are defined as follows:
		$$
		i^{p,q}(Y,w)=\sum _{\lambda \in \mathbb C}\sum _ki^{p,q+k}\mathbb H ^{p+q-1}(Y_\lambda ,
		\phi _{w-\lambda}\mathbb C _Y).
		$$
		
		\item
		In the general case denote by $j\colon Y\hookrightarrow Z$ the open embedding and define
		similarly
		$$
		i^{p,q}(Y,w)=\sum _{\lambda \in \mathbb C}\sum _ki^{p,q+k}\mathbb H ^{p+q-1}(Y_\lambda ,
		\phi _{w-\lambda}{\bf R}j_{*}\mathbb C _Y).
		$$
	\end{enumerate}
	Here $\mathbb{H}$ de notes the hypercohomology 
	of the constructible complexes of the
	sheaves of vanishing cycles, having higher derived images $\phi _{w-\lambda}{\bf R}j_{*}\mathbb C _Y$. 
\end{definition}

\subsection{The KKP conjecture}
The {\it KKP conjecture} states that the three invariants coincide, that is, 
$$f^{p,q}(Y,w)=h^{p,q}(Y,w)=i^{p,q}(Y,w).$$ 
For  $Y$  a specific rational surface with a map $w\colon Y\to \mathbb C$ such that the generic fibre is an elliptic curve 
\cite{LP} Lunts and Przyjalkowski proved the equality $f^{p,q}(Y,w)=h^{p,q}(Y,w)$ and gave an example  where~$i^{p,q}(Y,w)\neq h^{p,q}(Y,w)$. Thus, in full generality the conjecture is false. 
Nevertheless, 
Cheltsov and Przyjalkowski proved KKP the conjecture for Fano threefolds \cite{CP}
and Ballico, Gasparim, Rubilar, and San Martin proved  the KKP conjecture 
for minimal semisimple adjoint orbits \cite{BGRSM}. 
\begin{definition}In the cases when the KKP conjecture holds true, the invariants then define
	a new diamond, which we call the {\it KKP diamond}. 
\end{definition}

We now give some examples of Hodge diamonds and KKP diamonds computed in \cite{BGRSM}. 

\subsubsection{An example in 2  dimensions}
Consider  the semisimple adjoint orbit $\mathcal{O}_2$ of $\mathfrak{sl}(2, \mathbb C)$.
Hence, $\mathcal{O}_2$ can be viewed as  the affine hypersurface of $\mathbb C^3$ cut out by the equation
$x^2+yz=1$, see \cite[Sec.~2]{BBGGSM}.

The  Hodge diamond of $\mathcal{O}_2$ is:
\[ \begin{array}{ccccc}
& & 0 & \\   
& 0 &  & 0 & \\
\infty\!\! && 0 && 0  \\
&\infty\! & & 0  \\
&& \infty\!\! &
\end{array}\]
\[\,\quad\mathcal{O}_2 \]

We know by Cor. \ref{deforb} that $\mathcal O_2$ is a deformation of  $\mathrm{T}^*\mathbb{P}^1$.
The Hodge diamond of $\mathrm{T}^*\mathbb{P}^1$ is:
\[ \begin{array}{ccccc}
& & 0 & \\   
& 0 &  & 0 & \\
\infty\!\! && 1 && 0  \\
&\infty\! & & 0  \\
&& \infty\! &
\end{array}\]
\[\,\quad\mathrm{T}^*\mathbb{P}^1  \]
The KKP Hodge Diamond of  the Landau--Ginzburg model on $\mathrm{LG}(\mathcal O_2)$ obtained from Thm.~\ref{thm1} was calculated in  \cite[Sec.~7]{BGRSM} and is:
\[\begin{array}{cc}
\begin{array}{ccccc}
& & 0 & &\\   
& 0 &  & 0 & \\
0 && 2 && 0  \\
&0 & & 0 & \\
&& 0 & & \\
&&  & &
\end{array}\\
\mathrm{LG(2)}
\end{array}\]

\subsubsection{An example in 4  dimensions}
Our LG model is $(\mathcal{O}_3,f_H)=:\mathrm{LG}(3)$
obtained from Thm.~\ref{thm1} using the adjoint orbit  $\mathcal O_3$ of $\Diag(2,-1,-1)$
with potential $f_H$ corresponding to  the choice of $H= \Diag(1,0,-1)$.

The adjoint orbit $\mathcal O_3$ is a noncompact affine variety of dimension 4 that
has  the following Hodge diamond:
$$  \begin{array}{ccccccccc}
& & & & 0 & & \\
&  & & 0&  &0 & \\
& &0 & & 0 & &0 \\
&0&  &0  &  & 0&  &0 \\      
\infty\!\!&&0& & 0 &  & 0& & 0 \\
&\infty\! && 0 & & 0 && 0 \\
& &\infty\! & & 0 & & 0 \\
&&& \infty\! & & 0 & \\
&&& & \infty\! & &
\end{array}
$$
\vspace{-0.4cm}
\[\qquad\,\mathcal{O}_3\]

The  adjoint orbit ${\mathcal O_3}$ is a deformation of the cotangent bundle of 
$\mathbb{P}^2$, and $T^*\mathbb{P}^2$ has the following Hodge diamond:

$$  \begin{array}{rrccccccc}
& & & & 0 & & \\
&  & & 0&  &0 & \\
& &0 & & 0 & &0 \\
&0&  &0  &  & 0&  &0 \\      
\infty\!\!&&0& & 1 &  & 0& & 0 \\
&\infty\!\! && 0 & & 0 && 0 \\
& &\infty\!\! & & 1 & & 0 \\
&&& \infty\! & & 0 & \\
&&& & \infty\! & &
\end{array}
$$
\vspace{-0.3cm}

\[\quad\,\mathrm{T^*\mathbb P^2}\]

$ \mathrm{LG}(3) $ admits a  tame compactification, and 
the corresponding KKP diamond calculated in  \cite[Sec.~7]{BGRSM}  is:
$$  \begin{array}{rrrrcllll}
& & & &0& & & & \\
& & &0& &0& & & \\
& &0& &0& &0& & \\
&0& &0& &0& &0& \\      
0& &0& &3& &0& &0\\
&0& &0& &0& &0& \\
& &0& &0& &0& & \\
& & &0& &0& & & \\
& & & &0& & & &
\end{array}$$
\[\,\,\mathrm{LG}(3)\]

\subsubsection{The general case} \cite{BGRSM} proved the KKP conjecture for minimal adjoint orbits of 
$\mathfrak{sl}(n, \mathbb C)$. They calculated the 
KKP diamond of  $\mathrm{LG}(n)$ to be:

	$$\begin{array}{ccccccccc}

	& && &0 & && & \\
	& & & & & & & & \\
	& & &0 && 0& & & \\
	& & & &\vdots& & & & \\
	& &0 &\cdots &n&\cdots &0 & & \\
	& & & &\vdots& & & & \\
	& & &0 && 0&& & \\
	& & & & & & & & \\
	& && &0& && & \\

\end{array}$$

In conclusion, we have shown that our new concept of deformations of complex structures can be applied to many interesting examples.  Furthermore, we have described  classical Hodge theoretical invariants of cotangent bundles of projective spaces, and have compared them to the classical Hodge theoretical invariants of the nontrivial affine deformations of them, namely the minimal semisimple adjoint orbits of $\mathfrak{sl}(n, \mathbb C)$.
Finally we have described Landau--Ginzburg models on these adjoint orbits, and presented their KKP diamonds,
that is, diamonds containing the three new invariants defined by \cite{KKP}, which  as proved by \cite{BGRSM}, 
coincide for such orbits. How KKP diamonds vary under deformations is an interesting question that remains open. 
In fact, there are many delicate and intricate open questions about the deformation theory of Landau--Ginzburg models
both in complex and in symplectic geometry. 

\section{Acknowledgements}  We thank Severin Barmeier and the referee for suggesting several improvements to the text.

\bibliographystyle{amsalpha}

\end{document}